\documentclass{amsart}
\usepackage{amsfonts,amsmath,amsthm,amssymb}
\newtheorem{theorem}{Theorem}
\newtheorem{conjecture}{Conjecture}
\newtheorem{corollary}{Corollary}
\newtheorem{proposition}{Proposition}
\title[Integral group ring of the Mathieu  simple group $M_{24}$]
{Integral group ring of the \\ Mathieu  simple group $M_{24}$}
\author{V.A.~Bovdi,  A.B.~Konovalov}
\date{}

\address{
V.A.~Bovdi\newline
Institute of Mathematics, University of Debrecen \newline 
P.O.  Box 12, H-4010 Debrecen, Hungary \newline 
Institute of Mathematics and Informatics, 
College of Ny\'\i regyh\'aza \newline 
S\'ost\'oi \'ut 31/b, H-4410 Ny\'\i regyh\'aza, Hungary}
\email{vbovdi@math.klte.hu}

\address{
A.B.~Konovalov
\newline School of Computer Science, University of St Andrews,
\newline Jack Cole Building, North Haugh, St Andrews, Fife, KY16 9SX, Scotland}
\email{konovalov@member.ams.org}


\thanks{The research was supported by OTKA grants
No.T 037202, No.T 038059 and Francqui Stichting (Belgium) grant ADSI107}

\begin{document}
\begin{abstract}
We consider the  Zassenhaus conjecture  for the
normalized unit group of the integral group ring of  the
Mathieu sporadic group $M_{24}$. As a consequence, for this group we
confirm Kimmerle's conjecture on prime graphs.
\end{abstract}

\dedicatory{Dedicated to 70th birthday of Professor L.A.~Bokut}

\maketitle

\section{Introduction, conjectures   and main results}
\label{Intro}

Let $V(\mathbb Z G)$ be  the normalized units group of the integral
group ring $\mathbb Z G$ of  a finite group $G$. A famous Zassenhaus
conjecture \cite{Zassenhaus} says that every torsion unit $u\in V(\mathbb ZG)$ is
conjugate within the rational group algebra $\mathbb Q G$ to an
element in $G$.

For finite simple groups, the main tool for the investigation of the
Zassenhaus conjecture is the Luthar--Passi method, introduced in
\cite{Luthar-Passi} to solve this conjecture for $A_{5}$. Later
M.~Hertweck improved this method in \cite{Hertweck1} and used it for
the investigation of $PSL(2,F_{p^{n}})$. The Luthar--Passi method
proved to be useful for groups containing non-trivial normal
subgroups as well. Also some related properties
and  some weakened variations of the Zassenhaus  conjecture as well
can be found in \cite{Artamonov-Bovdi,Luthar-Trama} and
\cite{Bleher-Kimmerle,Kimmerle}. For some recent results we refer to
\cite{Bovdi-Hofert-Kimmerle,Bovdi-Konovalov,Hertweck2, Hertweck1,
Hertweck3, Hofert-Kimmerle}.

First of all, we need to introduce some notation. By $\# (G)$ we
denote the set of all primes dividing the order of $G$. The
Gruenberg--Kegel graph (or the prime graph) of $G$ is the graph $\pi
(G)$ with vertices labeled by the primes in $\# (G)$ and with an
edge from $p$ to $q$ if there is an element of order $pq$ in the
group $G$. The following weakened variation of the Zassenhaus 
conjecture was proposed in \cite{Kimmerle}:

\begin{conjecture}
{\bf (KC)} If $G$ is a
finite group then $\pi (G) =\pi (V(\mathbb Z G))$.
\end{conjecture}

In particular, in the same  paper  W.~Kimmerle verified that 
{\bf (KC)} holds for finite Frobenius and solvable groups. We remark 
that with respect to {\bf (ZC)} the investigation of Frobenius
groups was completed by  M.~Hertweck and the first author  in
\cite{Bovdi-Hertweck}. In \cite{Bovdi-Jespers-Konovalov,
Bovdi-Konovalov, Bovdi-Konovalov-M23, Bovdi-Konovalov-Linton,
Bovdi-Konovalov-Siciliano} {\bf (KC)} was confirmed
for the Mathieu simple groups $M_{11}$, $M_{12}$, $M_{22}$, $M_{23}$  and  the
sporadic Janko simple groups $J_1$, $J_{2}$ and $J_3$.

Here  we continue these investigations for the Mathieu simple group $M_{24}$. 
Despite using the Luthar--Passi method we are able to prove the rationally 
conjugacy only for torsion units of order 23 in $V(\mathbb Z M_{24})$, our main 
result gives a lot of information on partial augmentations of possible
torsion units and allows us to confirm {\bf (KC)} for the sporadic group $M_{24}$. 

It is well-known that the collection of conjugacy classes of $M_{24}$ is
\[
\begin{split}
\mathcal{C} =\{\; C_{1}, \; C_{2a},\; & C_{2b},\; C_{3a},\; C_{3b},\; C_{4a},\; C_{4b},\; C_{4c},\;\\ C_{5a},\; & C_{6a},\; C_{6b},\; C_{7a},\; C_{7b},\; C_{8a},\; C_{10a},\;  C_{11a},\; C_{12a}, \\
& \; C_{12b}, C_{14a},\; C_{14b},\; C_{15a},\; C_{15b},\; C_{21a},\; C_{21b},\; C_{23a},\; C_{23b} \; \},
\end{split}
\]
where the first
index denotes the order of the elements of this conjugacy class
and $C_{1}=\{ 1\}$. 
Suppose $u=\sum \alpha_g g \in V(\mathbb Z G)$
has finite order $k$. Denote by
$\nu_{nt}=\nu_{nt}(u)=\varepsilon_{C_{nt}}(u)=\sum_{g\in C_{nt}}
\alpha_{g}$, the partial augmentation of $u$ with respect to
$C_{nt}$. From the Berman--Higman Theorem 
(see \cite{Berman} and \cite{Sandling}, Ch.5, p.102) one knows that
$\nu_1 =\alpha_{1}=0$ and
\begin{equation}\label{E:1}
\sum_{C_{nt}\in \mathcal{C}} \nu_{nt}=1.
\end{equation}
Hence, for any character
$\chi$ of $G$, we get that $\chi(u)=\sum
\nu_{nt}\chi(h_{nt})$, where $h_{nt}$ is a
representative of a conjugacy class $ C_{nt}$.

The main result is the following.

\begin{theorem}\label{T:1}
Let $G$ denote the Mathieu  simple group $M_{24}$.  Let  $u$ be a torsion unit of
$V(\mathbb ZG)$ of order $|u|$ and let
\[
\begin{split}
\mathfrak{P}(u)=( \nu_{2a},\; &  \nu_{2b},\;  \nu_{3a},\;   \nu_{3b},\;   \nu_{4a},\;   \nu_{4b},\;
\nu_{4c},\;   \nu_{5a},\;  \nu_{6a},\;   \nu_{6b},\;   \nu_{7a},\;   \nu_{7b},\;   \nu_{8a},\;   \nu_{10a},\;  \\
&
\nu_{11a},\;   \nu_{12a},\;   \nu_{12b},\;   \nu_{14a},\;   \nu_{14b},\;   \nu_{15a},\;   \nu_{15b},\;   \nu_{21a},\;   \nu_{21b},\;
\nu_{23a},\;   \nu_{23b}\;) \in \mathbb Z^{25}
\end{split}
\]
be the tuple of partial augmentations of $u$.
The following properties hold.

\begin{itemize}

\item[(i)]
There is no elements of orders 
$22$, $33$, $35$, $46$, $55$, $69$, $77$, $115$, $161$ and $253$ in $V(\mathbb ZG)$.
Equivalently, if $|u| \not \in \{20, 24, 28, 30, 40, 42, 56, 60, 84, 120, 168 \}$,
then $|u|$ coincides  with the order of some element $g \in G$.
 
\item[(ii)]
If $|u| \in \{5,11, 23\}$, then $u$ is rationally
conjugate to some $g\in G$.

\item[(iii)]
If $|u|=2$, the tuple of the  partial augmentations of $u$ belongs
to the set
\[
\begin{split}
\big\{\;
( \mathfrak{P}(u) \in \mathbb Z^{25} \mid
    (\nu_{2a},\nu_{2b}) \in \{ \; ( 0, 1 ),\;   ( -2, 3 ),\;   ( 2, -1 ),\;   ( 1, 0 ), \; & \\
    ( 3, -2),\;   (-1, 2) \}, \quad \nu_{kx}=0,\; kx \not \in\{2a,2b\} \; & \big\} .
\end{split}
\]

\item[(iv)]
If $|u|=3$, the tuple of the  partial augmentations of $u$ belongs
to the set
\[
\begin{split}
\big\{\;
( \mathfrak{P}(u) \in \mathbb Z^{25} \mid
(\nu_{3a},\nu_{3b}) \in \{ \;
 (0, 1),\;   (2, -1),\;   (1, 0),\;   (3, -2),\;  &\\
  (-1, 2),\;   (4, -3) \; \},\quad \nu_{kx}=0,\; kx\not\in\{3a,3b\}\;& \big\} .
\end{split}
\]

\item[(v)]
If $|u|=7$, the tuple of the  partial augmentations of $u$ belongs
to the set
\[
\begin{split}
\big\{\;
( \mathfrak{P}(u) \in \mathbb Z^{25} \mid
(\nu_{7a},\nu_{7b}) \in \{ \; ( 0, 1),\; (2, -1),\; (1, 0),\; (-1, 2) \; \},\;&\\
\nu_{kx}=0,\;kx\not\in\{7a,7b\}\;&\big\}.
\end{split}
\]

\item[(vi)]
If $|u|=10$, the tuple of the  partial augmentations of $u$ belongs
to the set
\[
\begin{split}
\big\{ \quad ( \mathfrak{P}(u) \in \mathbb Z^{25} \mid
(\nu_{2a},\nu_{2b}, \nu_{5a},\nu_{10a}) \in \{ \; (-3, 1, 5, -2), (-2, 0, 5, -2), & \\
(-2, 2, 5, -4),(-1, -1, 5, -2), (-1, 1, 5, -4), (0, -2, 0, 3), (0, 0, 0, 1), & \\
(0, 2, 0, -1), (1, -1, 0, 1), (1, 1, 0, -1),(1, 3, 0, -3) \; \}, & \\
\nu_{kx}=0,\;kx\not\in\{\nu_{2a},\nu_{2b}, \nu_{5a},\nu_{10a}\}& \; \big\}.
\end{split}
\]

\end{itemize}
\end{theorem}

Note that using our implementation of the Luthar--Passi method,
which we intend to make available in the GAP package LAGUNA \cite{LAGUNA},
it is possible to compute 34 possible tuples of partial augmentations
for units of order 15 and 21 tuple for units of order 21 listed in the Appendix.

As an immediate consequence of the  part (i) of the Theorem we
obtain

\begin{corollary} If $G=M_{24}$ then
$\pi(G)=\pi(V(\mathbb ZG))$.
\end{corollary}

\section{Preliminaries}

The following result relates the solution of
the Zassenhaus conjecture to partial augmentations
of torsion units.

\begin{proposition}\label{P:5}
(see \cite{Luthar-Passi} and
Theorem 2.5 in \cite{Marciniak-Ritter-Sehgal-Weiss})
Let $u\in V(\mathbb Z G)$
be of order $k$. Then $u$ is conjugate in $\mathbb
QG$ to an element $g \in G$ if and only if for
each $d$ dividing $k$ there is precisely one
conjugacy class $C$ with partial augmentation
$\varepsilon_{C}(u^d) \neq 0 $.
\end{proposition}

The next result already yield that several
partial augmentations are zero.

\begin{proposition}\label{P:4}
(see \cite{Hertweck2}, Proposition 3.1;
\cite{Hertweck1}, Proposition 2.2)
Let $G$ be a finite
group and let $u$ be a torsion unit in $V(\mathbb
ZG)$. If $x$ is an element of $G$ whose $p$-part,
for some prime $p$, has order strictly greater
than the order of the $p$-part of $u$, then
$\varepsilon_x(u)=0$.
\end{proposition}

The key restriction on partial augmentations is given 
by the following result that is the cornerstone of
the Luthar--Passi method.

\begin{proposition}\label{P:1}
(see \cite{Hertweck1, Luthar-Passi})
Let either $p=0$ or $p$ a prime divisor of $|G|$.
Suppose that $u\in V( \mathbb Z G) $ has finite
order $k$ and assume $k$ and $p$ are coprime in
case $p\neq 0$. If $z$ is a complex primitive $k$-th root
of unity and $\chi$ is either a classical
character or a $p$-Brauer character of $G$ then,
for every integer $l$, the number
\begin{equation}\label{E:2}
\mu_l(u,\chi, p )=\textstyle\frac{1}{k} \sum_{d|k}Tr_{ \mathbb Q (z^d)/ \mathbb Q }
\{\chi(u^d)z^{-dl}\}
\end{equation}
is a non-negative integer.
\end{proposition}

Note that if $p=0$, we will use the notation $\mu_l(u,\chi
, * )$ for $\mu_l(u,\chi , 0)$.

Finally, we shall use the well-known bound for
orders of torsion units.

\begin{proposition}\label{P:2}  (see  \cite{Cohn-Livingstone})
The order of a torsion element $u\in V(\mathbb ZG)$
is a divisor of the exponent of $G$.
\end{proposition}

\section{Proof of the Theorem}

Throughout this section we denote $M_{24}$ by $G$. It is well known
\cite{GAP} that $|G|=2^{10}\cdot 3^{3}\cdot 5 \cdot 7
\cdot 11\cdot 23$ and $exp(G)=2^{3}\cdot 3\cdot 5 \cdot 7 \cdot
11\cdot 23$. The character table of $G$, as well as the $p$-Brauer
character tables, where $p \in \{2, 3, 5, 7, 11, 23\}$, 
can be found using the computational algebra system GAP \cite{GAP}, 
which derives these data from \cite{AFG,ABC}. 
Throughout the paper we will use the notation, inclusive
the indexation, for the characters and conjugacy classes as used in
the GAP Character Table Library.

Since the group $G$ possesses elements of orders $2$, $3$, $4$, $5$,
$6$, $7$, $8$, $10$, $11$, $12$, $14$, $15$, $21$ and $23$,
first of all we will investigate  units of some of these orders (except
units of orders $4$, $6$, $8$, $12$ and $14$). After this, by
Proposition \ref{P:2}, the order of each torsion unit divides the
exponent of $G$, and in the first instance we should consider 
units of orders $20$, $22$, $24$, $28$, $30$, $33$, $35$, $42$,  
$46$, $55$, $69$, $77$, $115$, $161$ and $253$. We will omit orders
$20$, $24$, $28$, $30$ and $42$ that do not contribute to {\bf (KC)}, 
and this enforces us to add to the list of exceptions in part (i) of Theorem also orders
$40$, $56$, $60$, $84$, $120$ and $168$, but no more because
of restrictions imposed by the exponent of $G$.

Thus, we will prove that units of orders 
$22$, $33$, $35$, $46$, $55$, $69$, $77$, $115$, $161$ and $253$
do not appear in $V(\mathbb ZG)$.

Now we consider each case separately.

\noindent$\bullet$ Let $u$ be an involution. By (\ref{E:1}) and
Proposition \ref{P:4} we have that $\nu_{2a}+\nu_{2b}=1$. Applying
Proposition \ref{P:1} to characters $\chi_{2}$ we get the following system
\[
\begin{split}
\mu_{0}(u,\chi_{2},*) & = \textstyle \frac{1}{2} (7 \nu_{2a} -  \nu_{2b} +
23) \geq 0; \quad
\mu_{1}(u,\chi_{2},*)   = \textstyle \frac{1}{2} (-7 \nu_{2a} +  \nu_{2b} +
23) \geq 0. 
\end{split}
\]
From these restrictions and the requirement that all
$\mu_i(u,\chi_{j},p)$ must be non-ne\-ga\-tive integers
we get the six  pairs $(\nu_{2a},\nu_{2b})$ listed in
part (iii) of our Theorem.

\noindent $\bullet$ Let $u$ be a unit of order $3$. By (\ref{E:1})
and Proposition \ref{P:4} we get $\nu_{3a}+\nu_{3b}=1$. By (\ref{E:2}) we obtain  the system of inequalities
\[
\begin{split}
\mu_{0}(u,\chi_{2},*) & = \textstyle \frac{1}{3} (10 \nu_{3a} - 2 \nu_{3b} +
23) \geq 0;\\ 
\mu_{1}(u,\chi_{2},*) & = \textstyle \frac{1}{3} (-5 \nu_{3a} +  \nu_{3b} +
23) \geq 0. 
\end{split}
\]
Clearly, using  the condition for
$\mu_i(u,\chi_{j},p)$ to be non-negative integers,
we obtain six pairs $(\nu_{3a},\nu_{3b})$ listed in
part (iv) of the Theorem \ref{T:1}.

\noindent $\bullet$ Let $u$ be a unit of order either $5$ or $11$.
Using Proposition \ref{P:4} and (\ref{E:2}) we obtain that all partial
augmentations except one are zero. Thus by Proposition \ref{P:4} the particular proof of
part (ii) of the Theorem \ref{T:1} is done.

\noindent $\bullet$ Let $u$ be a unit of order $7$. By (\ref{E:1})
and Proposition \ref{P:4} we get $\nu_{7a}+\nu_{7b}=1$. By (\ref{E:2}) we obtain  the system of inequalities
\[
\begin{split}
\mu_{1}(u,\chi_{3},*) & = \textstyle \frac{1}{7} (4 \nu_{7a} - 3 \nu_{7b} +
45) \geq 0;\\ 
\mu_{1}(u,\chi_{2},2) & = \textstyle \frac{1}{7} (-4 \nu_{7a} + 3 \nu_{7b} +
11) \geq 0. 
\end{split}
\]
Again, using  the condition for $\mu_i(u,\chi_{j},p)$ to be
non-negative integers, we obtain four pairs $(\nu_{7a},\nu_{7b})$
listed in part (v) of the Theorem \ref{T:1}.

\noindent$\bullet$ Let $u$ be a unit of order $10$. By (\ref{E:1})
and Proposition \ref{P:4} we have that
\begin{equation}\label{E:3}
\nu_{2a}+\nu_{2b}+\nu_{5a}+\nu_{10a}=1.
\end{equation}
Since $|u^5|=2$, we need to consider six cases defined by part (iii) of the Theorem \ref{T:1}.

\noindent
Case 1. Let $\chi(u^5)=\chi(2a)$. Put
\begin{equation}\label{E:4}
\begin{split}
t_1 & = 7 \nu_{2a} -  \nu_{2b} + 3 \nu_{5a} -  \nu_{10a}, 
\qquad 
t_2 = 3 \nu_{2a} - 5 \nu_{2b}, \\
t_3 & = 14 \nu_{2a} + 6 \nu_{2b} +  \nu_{5a} +  \nu_{10a}.
\end{split}
\end{equation}
Applying Proposition \ref{P:1}, we get the system with indeterminates $t_1$, $t_2$ and $t_3$
\[
\begin{split}
\mu_{0}(u,\chi_{2},*) & = \textstyle \frac{1}{10} (4 t_1 + 42) \geq 0; \quad 
\mu_{5}(u,\chi_{2},*)   = \textstyle \frac{1}{10} (-4 t_1 + 28) \geq 0;\\ 
\mu_{5}(u,\chi_{3},*) & = \textstyle \frac{1}{10} (4 t_2 + 48) \geq 0; \quad 
\mu_{0}(u,\chi_{3},*)   = \textstyle \frac{1}{10} (-12 t_2 + 42) \geq 0; \\
\mu_{0}(u,\chi_{7},*) & = \textstyle \frac{1}{10} (8 t_3 + 288) \geq 0; \quad 
\mu_{5}(u,\chi_{7},*)   = \textstyle \frac{1}{10} (-8 t_3 + 232) \geq 0.\\ 
\end{split}
\]
Its solution are $t_1 \in \{-8, -3, 2, 7\}$,
$t_2 \in \{-12, -7, -2, 3, 8\}$
and
$$
t_3 \in \{ -36, -31, -26, -21, -16, -11, -6, -1, 4, 9, 14, 19, 24, 29 \}.
$$
Substituting values of $t_1$, $t_2$ and $t_3$ in (\ref{E:4}), and adding the condition (\ref{E:3}),
we obtain the system of linear equations for $\nu_{2a}$, $\nu_{2b}$, $\nu_{5a}$, and $\nu_{10a}$.
Since 
$ \tiny { \left| 
\begin{matrix} 
1&1&1&1\\
7&-1&3&-1\\
3&-5&0&0\\
14&6&1&1\\
\end{matrix} 
\right| }  
\not=0$, this system has the unique solution for each $t_1$, $t_2$, $t_3$, and the only 
integer solutions are $(1, -1,0, 1)$,\quad  $(1, 1, 0, -1)$ and $(1, 3, 0, -3)$.

\noindent
Case 2. Let $\chi(u^5)=\chi(2b)$. Put  $t_1 = 7 \nu_{2a} -  \nu_{2b} + 3 \nu_{5a} -  \nu_{10a}$, $t_2 = 3 \nu_{2a} - 5 \nu_{2b}$ and
$t_3 = 14 \nu_{2a} + 6 \nu_{2b} +  \nu_{5a} +  \nu_{10a}$.
Again using Proposition \ref{P:1}, we obtain that
\[
\begin{split}
\mu_{0}(u,\chi_{2},*) & = \textstyle \frac{1}{10} (4 t_1 + 34) \geq 0; \quad 
\mu_{5}(u,\chi_{2},*)   = \textstyle \frac{1}{10} (-4 t_1 + 36) \geq 0;\\ 
\mu_{5}(u,\chi_{3},*) & = \textstyle \frac{1}{10} (4 t_2 + 40) \geq 0; \quad 
\mu_{0}(u,\chi_{3},*)   = \textstyle \frac{1}{10} (-4 t_2 + 50) \geq 0; \\
\mu_{0}(u,\chi_{7},*) & = \textstyle \frac{1}{10} (8 t_3 + 272) \geq 0; \quad 
\mu_{5}(u,\chi_{7},*)   = \textstyle \frac{1}{10} (-8 t_3 + 248) \geq 0.\\ 
\end{split}
\]
From this follows that
$t_1 \in \{-6, -1, 4, 9 \}$,
$t_2 \in \{-10, -5, 0, 5, 10 \}$
and
$$
t_3 \in \{-34, -29, -24, -19, -14, -9, -4, 1, 6, 11, 16, 21, 26, 31 \}.
$$
Using the same considerations as in the previous case, we obtain
only three solutions
$(0, -2, 0, 3)$, $( 0, 0, 0, 1)$ and $( 0, 2, 0, -1)$ that satisfy
these restrictions and the condition that $\mu_i(u,\chi_{j},p)$ 
are non-negative integers.

\noindent
Case 3. Let $\chi(u^5)=-2\chi(2a)+3\chi(2b)$.
Put  
$t_1 = 7 \nu_{2a} -  \nu_{2b} + 3 \nu_{5a} -  \nu_{10a}$, 
$t_2 = 3 \nu_{2a} - 5 \nu_{2b}$
and 
$t_3 = 14 \nu_{2a} + 6 \nu_{2b} +  \nu_{5a} +  \nu_{10a}$.
As before, by Proposition \ref{P:1}, we obtain that
\[
\begin{split}
\mu_{0}(u,\chi_{2},*) & = \textstyle \frac{1}{10} (4 t_1 + 18) \geq 0; \quad 
\mu_{2}(u,\chi_{2},*)   = \textstyle \frac{1}{10} (- t_1 + 3) \geq 0; \\ 
\mu_{5}(u,\chi_{3},*) & = \textstyle \frac{1}{10} (4 t_2 + 24) \geq 0; \quad 
\mu_{0}(u,\chi_{3},*)   = \textstyle \frac{1}{10} (-t t_2 + 66) \geq 0; \\ 
\mu_{0}(u,\chi_{7},*) & = \textstyle \frac{1}{10} (8 t_3 + 240) \geq 0; \quad 
\mu_{5}(u,\chi_{7},*)   = \textstyle \frac{1}{10} (-8 t_3 + 280) \geq 0. \\
\end{split}
\]
From the last system of inequalities, we get
$t_1=3$,
$t_2 \in \{-6, -1, 4, 9, 14 \}$
and
$$
t_3 \in \{ -30, -25, -20, -15, -10, -5, 0, 5, 10, 15, 20, 25, 30, 35 \},
$$
and using the same considerations as in the previous case, 
we deduce that there is only one solution $(-2, 0, 5, -2)$ 
satisfying the previous restrictions and the condition that
$\mu_i(u,\chi_{j},p)$ are non-negative integers.

\noindent
Case 4. Let $\chi(u^5)=2\chi(2a)-\chi(2b)$. Again, for the same
$t_1$, $t_2$ and $t_3$ we have
\[
\begin{split}
\mu_{1}(u,\chi_{2},*) & = \textstyle \frac{1}{10} (t_1 + 5) \geq 0;\qquad \; 
\mu_{5}(u,\chi_{2},*)   = \textstyle \frac{1}{10} (-4 t_1 + 20) \geq 0;\\ 
\mu_{5}(u,\chi_{3},*) & = \textstyle \frac{1}{10} (4 t_2 + 56) \geq 0; \quad \; 
\mu_{0}(u,\chi_{3},*)   = \textstyle \frac{1}{10} (-4 t_2 + 34) \geq 0;\\ 
\mu_{0}(u,\chi_{7},*) & = \textstyle \frac{1}{10} (8 t_3 + 304) \geq 0; \quad 
\mu_{5}(u,\chi_{7},*)   = \textstyle \frac{1}{10} (-8 t_3 + 216) \geq 0. 
\end{split}
\]
It follows that \quad $t_1\in\{-5,5\}$, \quad $t_2 \in \{-14, -9, -4, 1, 6\}$ \quad
and
$$
t_3 \in \{ -38, -33, -28, -23, -18, -13, -8, -3, 2, 7, 12, 17, 22, 27 \},
$$
and we obtain tree solutions 
$\{ \; (-3, 0, 5, -1),\; (-3, 1, 5,-2),\; (2, 0, -5, 4) \; \}$ satisfying 
the inequalities above. Now using the following additional inequalities:
\[
\begin{split}
\mu_{1}(u,\chi_{3},*) & = \textstyle \frac{1}{10} (-3 \nu_{2a} + 5 \nu_{2b} +
56) \geq 0; \\
\mu_{5}(u,\chi_{5},11) & = \textstyle \frac{1}{10} (-84 \nu_{2a} -
52 \nu_{2b} + 4 \nu_{5a} - 12 \nu_{10a} + 196) \geq 0, \\
\end{split}
\]
it remains only one solution $( -3, 1, 5, -2)$.

\noindent
Case 5. Let $\chi(u^5)=3\chi(2a)-2\chi(2b)$. \quad Put \quad
$t_1 = 7 \nu_{2a} -  \nu_{2b} + 3 \nu_{5a} -  \nu_{10a}$, \qquad
$t_2 = 3 \nu_{2a} - 5 \nu_{2b}$
\; and \;
$t_3 = 14 \nu_{2a} + 6 \nu_{2b} +  \nu_{5a} +  \nu_{10a}$.
Again by (\ref{E:2}) we obtain that
\[
\begin{split}
\mu_{1}(u,\chi_{2},*) & = \textstyle \frac{1}{10} (t_1 - 3) \geq 0; \qquad 
\mu_{5}(u,\chi_{2},*)   = \textstyle \frac{1}{10} (-4 t_1 + 12) \geq 0; \\ 
\mu_{5}(u,\chi_{3},*) & = \textstyle \frac{1}{10} (4 t_2 + 64) \geq 0; \quad 
\mu_{0}(u,\chi_{3},*)   = \textstyle \frac{1}{10} (-4 t_2 + 26) \geq 0; \\ 
\mu_{0}(u,\chi_{7},*) & = \textstyle \frac{1}{10} (8 t_3 + 320) \geq 0; \quad
\mu_{5}(u,\chi_{7},*)   = \textstyle \frac{1}{10} (-8 t_3 + 200) \geq 0. 
\end{split}
\]
It is easy to check that $t_1=3$, $t_2\in \{-16, -11, -6, -1, 4\}$
and
$$
t_3 \in \{  -40, -35, -30, -25, -20, -15, -10, -5, 0, 5, 10, 15, 20, 25 \}.
$$
So we obtained the following five solutions:
$$\{ (-2, -2, 5, 0), (-2, -1, 5, -1),
(-2, 0, 5, -2), (-2, 1, 5, -3), (-2, 2, 5, -4)\}.
$$
Now after using the following two additional inequalities:
\[
\begin{split}
\mu_{1}(u,\chi_{3},*) & = \textstyle \frac{1}{10} (-3 \nu_{2a} + 5 \nu_{2b} +
64) \geq 0;\\ 
\mu_{0}(u,\chi_{5},11) & = \textstyle \frac{1}{10} (84 \nu_{2a} +
52 \nu_{2b} - 4 \nu_{5a} + 12 \nu_{10a} + 262) \geq 0,\\ 
\end{split}
\]
it remains only two  solutions 
$\{ \; ( -2, 0, 5, -2),\; (-2, 2, 5, -4) \; \}$.

\noindent
Case 6. Let $\chi(u^5)=-\chi(2a)+2\chi(2b)$.
Put \quad $t_1 = 7 \nu_{2a} -  \nu_{2b} + 3 \nu_{5a} -  \nu_{10a}$, \qquad 
$t_2 = 3 \nu_{2a} - 5 \nu_{2b}$ and
$t_3 = 14 \nu_{2a} + 6 \nu_{2b} +  \nu_{5a} +  \nu_{10a}$.
Similarly, we get
\[
\begin{split}
\mu_{0}(u,\chi_{2},*) & = \textstyle \frac{1}{10} ( 4 t_1 + 26) \geq 0; \quad 
\mu_{2}(u,\chi_{2},*)   = \textstyle \frac{1}{10} (- t_1 + 11) \geq 0; \\ 
\mu_{5}(u,\chi_{3},*) & = \textstyle \frac{1}{10} (4 t_2 + 32) \geq 0; \quad 
\mu_{0}(u,\chi_{3},*)   = \textstyle \frac{1}{10} (-4 t_2 + 58) \geq 0; \\ 
\mu_{0}(u,\chi_{7},*) & = \textstyle \frac{1}{10} (8 t_3 + 256) \geq 0; \quad 
\mu_{5}(u,\chi_{7},*)   = \textstyle \frac{1}{10} (-8 t_3 + 264) \geq 0. 
\end{split}
\]
We have the following restrictions:
$t_1\in \{1, 11\}$,
$t_2 \in \{-8, -3, 2, 7, 12\}$
and
$$
t_3 \in \{ -32, -27, -22, -17, -12, -7, -2, 3, 8, 13, 18, 23, 28, 33 \},
$$
that lead to the following five solutions
$$
\{(-1, -3, 5, 0 ),\; ( -1, -2, 5, -1 ),\;  (-1, -1, 5, -2 ),\;  (-1, 0, 5, -3 ),\;  ( -1, 1, 5, -4)\}
$$
which satisfy the above inequalities. After considering two additional inequalities
\[
\begin{split}
\mu_{1}(u,\chi_{3},*) & = \textstyle \frac{1}{10} (-3 \nu_{2a} + 5 \nu_{2b} +
32) \geq 0;\\ 
\mu_{0}(u,\chi_{5},11) & = \textstyle \frac{1}{10} (84 \nu_{2a} +
52 \nu_{2b} - 4 \nu_{5a} + 12 \nu_{10a} + 230) \geq 0,\\ 
\end{split}
\]
only two solutions remains: $\{ \; ( -1, -1, 5, -2), \; (-1, 1, 5, -4) \; \}$.

Thus, the union of solutions for all six cases gives us
part (vi) of the Theorem.

\noindent$\bullet$ Let $u$ be a unit of order $15$. 
By (\ref{E:1})
and Proposition \ref{P:4} we obtain that
$$
\nu_{3a}+\nu_{3b}+\nu_{5a}+\nu_{15a}+\nu_{15b}=1.
$$
Since $|u^5|=3$, according to part (iv) of the Theorem
we need to consider six cases. Using the LAGUNA package
\cite{LAGUNA}, in all of them we constructed and solved 
systems of inequalities that give us 34 solutions listed in 
the Appendix.

\noindent$\bullet$ Let $u$ be a unit of order $21$. By (\ref{E:1})
and Proposition \ref{P:4} we obtain that
$$
\nu_{3a}+\nu_{3b}+\nu_{7a}+\nu_{7b}+\nu_{21a}+\nu_{21b}=1.
$$
We need to consider $24$ cases determined by parts (iv) and (v) of the Theorem \ref{T:1}.
We write down explicitly  the details of the first case, the treatment of the other ones
are similar. Our computation was helped by the LAGUNA package \cite{LAGUNA}.

Let  $\chi(u^3)=\chi(7a)$ and $\chi(u^7)=\chi(3a)$, for any
character $\chi$ of $G$. Put
\[
\begin{split}
t_1&= 5 \nu_{3a} -  \nu_{3b} + 2 \nu_{7a} + 2 \nu_{7b} - \nu_{21a}
- \nu_{21b}, \qquad \\
t_2 &=6 \nu_{3b} -  \nu_{7a} -  \nu_{7b} - \nu_{21a} -
\nu_{21b},\quad  t_3 = 3 \nu_{3b} + 3 \nu_{7a} - 4 \nu_{7b} + 3
\nu_{21a} - 4 \nu_{21b},\\
t_4 & =  \nu_{3a}, \quad 
t_5   = 3 \nu_{3b} - 6 \nu_{7a} + 8 \nu_{7b} + 3 \nu_{21a} - 4 \nu_{21b}.
\end{split}
\]
Applying Proposition \ref{P:1} to characters $\chi_{2}$, $\chi_{3}$,
$\chi_{4}$, $\chi_{7}$ and $\chi_{15}$ we get
\[
\begin{split}
\mu_{0}(u,\chi_{2},*) & = \textstyle \frac{1}{21} (5 t_1 + 45) \geq 0;\qquad 
\mu_{7}(u,\chi_{2},*)   = \textstyle \frac{1}{21} (-6 t_1 + 30) \geq 0;\\ 
\mu_{0}(u,\chi_{3},*) & = \textstyle \frac{1}{21} (6 t_2 + 42) \geq 0; \qquad
\mu_{7}(u,\chi_{3},*)   = \textstyle \frac{1}{21} (-3 t_2 + 42) \geq 0;\\ 
\mu_{1}(u,\chi_{3},*) & = \textstyle \frac{1}{21} (t_3 + 49) \geq 0; \;\qquad 
\mu_{9}(u,\chi_{3},*)   = \textstyle \frac{1}{21} (-2 t_3 + 49) \geq 0;\\ 
\mu_{0}(u,\chi_{7},*) & = \textstyle \frac{1}{21} (108 t_4 + 270) \geq 0;\ \; 
\mu_{7}(u,\chi_{7},*)   = \textstyle \frac{1}{21} (-54 t_4 + 243) \geq 0;\\ 
\mu_{9}(u,\chi_{15},*) & = \textstyle \frac{1}{21} ( 2 t_5 + 1043) \geq 0; \quad 
\mu_{1}(u,\chi_{15},*)   = \textstyle \frac{1}{21} (-t_5 + 1043) \geq 0. \\
\end{split}
\]
Solution of this system of inequalities gives \quad
$t_1 \in \{ -2, 5 \}$, \quad $t_2 \in \{-7, 0, 7, 14\}$, \quad $t_3 \in \{ -49, -28, -7, 14 \}$,\quad  $t_4=1$
\quad and \quad
$t_5 \in \{ 14 + 21k \mid -25 \leq k \leq 49 \}$.

Using computer we get $1200$ solutions satisfying inequalities above.

After considering the following four additional inequalities
\[
\begin{split}
\mu_{9}(u,\chi_{2},2) & = \textstyle \frac{1}{21} (-4 \nu_{3a} + 2 \nu_{3b} +
6 \nu_{7a} - 8 \nu_{7b} - 12 \nu_{21a} + 16 \nu_{21b} + 11) \geq 0;\\ 
\mu_{1}(u,\chi_{2},2) & = \textstyle \frac{1}{21} (2 \nu_{3a} -  \nu_{3b} -
3 \nu_{7a} + 4 \nu_{7b} + 6 \nu_{21a} - 8 \nu_{21b} + 5) \geq 0;\\ 
\mu_{0}(u,\chi_{4},2) & = \textstyle \frac{1}{21} (-12 \nu_{3a} +
24 \nu_{3b} - 18 \nu_{7a} - 18 \nu_{7b} - 18 \nu_{21a} - 18 \nu_{21b} +
33) \geq 0;\\ 
\mu_{3}(u,\chi_{3},*) & = \textstyle \frac{1}{21} (- 6 \nu_{3b} +
8 \nu_{7a} - 6 \nu_{7b} + 8 \nu_{21a} - 6 \nu_{21b} + 42) \geq 0,\\ 
\end{split}
\]
it remains only two  solutions: 
$\{ \; (1, 2, -1, 1, -2, 0),\; (1, 2, 1, -1, -1, -1) \; \}$.

Similarly, using the LAGUNA package \cite{LAGUNA} we can construct the system
of inequalities for the remaining 23 cases. The union of all solutions
give us the list of solutions given in the Appendix.

\noindent$\bullet$ Let $u$ be a unit of order $23$. By (\ref{E:1})
and Proposition \ref{P:4} we get $\nu_{23a}+\nu_{23b}=1$.
By (\ref{E:2}) we obtain the following system of inequalities
\[
\begin{split}
\mu_1(u,\chi_{10},*) & =\textstyle \frac{1}{23} ( 12 \nu_{23a} -11 \nu_{23b} + 770 ) \geq 0;\\
\mu_5(u,\chi_{10},*) & = \textstyle\frac{1}{23} ( -11 \nu_{23a} + 12 \nu_{23b} + 770 ) \geq 0;\\
\mu_1(u,\chi_2,2) & = \textstyle\frac{1}{23} ( 12 \nu_{23a} -11 \nu_{23b} + 11 ) \geq 0;\\
\mu_5(u,\chi_2,2) & = \textstyle\frac{1}{23} ( -11 \nu_{23a} + 12 \nu_{23b} + 11 ) \geq 0;\\
\mu_1(u,\chi_7,2) & = \textstyle\frac{1}{23} ( -13 \nu_{23a} + 10 \nu_{23b} + 220 ) \geq 0;\\
\mu_5(u,\chi_7,2) & = \textstyle\frac{1}{23} ( 10 \nu_{23a} -13 \nu_{23b} + 220 ) \geq 0;\\
\mu_1(u,\chi_{10},2) & = \textstyle\frac{1}{23} ( 25 \nu_{23a} -21 \nu_{23b} + 320 ) \geq 0;\\
\mu_5(u,\chi_{10},2) & = \textstyle\frac{1}{23} ( -21 \nu_{23a} + 25 \nu_{23b} + 320 ) \geq 0,\\
\end{split}
\]
which has only two trivial solutions $(\nu_{23a},\nu_{23b})\in\{(1,0), (0,1)\}$.
Thus, by Proposition \ref{P:5} we conclude that each torsion unit of order 23 
is rationally conjugate to some $g \in G$, and this completes the proof of 
part (ii) of the Theorem.

\noindent$\bullet$ Let $u$ be a unit of order $22$. By (\ref{E:1})
and Proposition \ref{P:4} we have that
$$
\nu_{2a}+\nu_{2b}+\nu_{11a}=1.
$$
Since $|u^{11}|=2$,  we need to
consider six cases for any character $\chi$ of $G$. They are defined by part (iii) of the Theorem. Put
\begin{equation}\label{E:5}
(\alpha,\beta,\gamma,\delta)= {\tiny \begin{cases}
(40,26,58,52), &\quad \text{if }\quad  \chi(u^{11})=\chi(2a);\\
(32,34,50,60), &\quad \text{if }\quad  \chi(u^{11})=\chi(2b); \\
(16,5,34,76),  &\quad \text{if }\quad  \chi(u^{11})=-2\chi(2a)+\chi(2b);\\
(48,18, 66,44),&\quad \text{if }\quad  \chi(u^{11})=2\chi(2a)-\chi(2b);\\
(-1, 10, 74, 36),&\quad \text{if }\quad  \chi(u^{11})=3\chi(2a)-2\chi(2b);\\
(24, 42, 42,68),&\quad \text{if }\quad  \chi(u^{11})=-\chi(2a)+2\chi(2b),\\
\end{cases}}
\end{equation}
\begin{equation}\label{E:6}
t_1 = 7 \nu_{2a} -  \nu_{2b} +  \nu_{11a} \qquad \text{and} \qquad  t_2 = 3 \nu_{2a} - 5 \nu_{2b} -  \nu_{11a}.
\end{equation}
If  $\chi(u^{11})=3\chi(2a)-2\chi(2b)$, by  (\ref{E:2}) we obtain the system
\begin{equation}\label{E:5}
\begin{split}
\mu_{0}(u,\chi_{2},*)  & = \textstyle \frac{1}{22} (10 t_1 + \alpha) \geq 0;\quad 
\mu_{11}(u,\chi_{2},*)   = \textstyle \frac{1}{22} (-10 t_1 + \beta) \geq 0;\\ 
\mu_{11}(u,\chi_{3},*) & = \textstyle \frac{1}{22} (10 t_2 + \gamma) \geq 0;\quad 
\mu_{0}(u,\chi_{3},*)    = \textstyle \frac{1}{22} (-10 t_2 + \delta) \geq 0.\\ 
\end{split}
\end{equation}
For each of the cases of (\ref{E:5}), we solve the system (\ref{E:6})
for $t_1$ and $t_2$. Then we obtain the following six solutions.
\begin{itemize}
\item[(i)] \quad  $\chi(u^{11})=\chi(2a)$. We get $t_1=4$ and $t_2=3$.
\item[(ii)] \quad  $\chi(u^{11})=\chi(2b)$. We get $t_1=-1$ and $t_2\in\{-5,6\}$.
We have the solution $(\nu_{2a},\nu_{2b}, \nu_{11a})=(0, 1, 0)$. After considering the additional
restriction $\mu_{1}(u,\chi_{5},*) = \textstyle \frac{1}{22} (7 \nu_{2a} - 9 \nu_{2b} +
240)=\textstyle \frac{231}{22}$. Since $\mu_{1}(u,\chi_{5},*)$ is not an integer, we obtain
a contradiction, so in this case there is no solution.
\item[(iii)] \quad  $\chi(u^{11})=-2\chi(2a)+\chi(2b)$. We get $t_1=5$ and $t_2=1$.
\item[(iv)] \quad  $\chi(u^{11})=2\chi(2a)-\chi(2b)$. In this case there is no solution for $t_1$.
\item[(v)] \quad $\chi(u^{11})=\chi(u^{11})=3\chi(2a)-2\chi(2b)$. We
get $t_1=1$ and $t_2=-3$.
\item[(vi)] \quad  $\chi(u^{11})=-\chi(2a)+2\chi(2b)$. We get $t_1=2$ and $t_2=-2$.
\end{itemize}
Finally, assume that  $\chi(u^{11})=3\chi(2a)-2\chi(2b)$. Put\quad  $t_1 = 7 \nu_{2a} -  \nu_{2b} +  \nu_{11a}$\quad
and $t_2 = 3 \nu_{2a} - 5 \nu_{2b} -  \nu_{11a}$.
Again, by  (\ref{E:2}) we obtain the system of inequalities
\[
\begin{split}
\mu_{1}(u,\chi_{2},*) & = \textstyle \frac{1}{22} ( t_1 - 1) \geq 0;\qquad \; 
\mu_{11}(u,\chi_{2},*)  = \textstyle \frac{1}{22} (-10 t_1 + 10) \geq 0;\\ 
\mu_{11}(u,\chi_{3},*)& = \textstyle \frac{1}{22} (10 t_2 + 74) \geq 0;\quad 
\mu_{0}(u,\chi_{3},*)   = \textstyle \frac{1}{22} (-10 t_2 + 36) \geq 0,\\ 
\end{split}
\]
with  integral solution  $(t_1,t_2)=(1 , -3)$. Now we substitute the obtained
values of $t_1$ and $t_2$ into the system of equations (\ref{E:4}). Then we can conclude
that it is impossible to find integer solution of (\ref{E:4}) for $\nu_{2a}$, $\nu_{2b}$ and $\nu_{11a}$.

\noindent$\bullet$ Let $u$ be a unit of order $33$. By (\ref{E:1})
and Proposition \ref{P:4} we have that
$$
\nu_{3a}+\nu_{3b}+\nu_{11a}=1.
$$
Since $|u^{11}|=3$, for any character $\chi$ of $G$ we need to
consider six cases, defined by part (iv) of the Theorem. Put
\begin{equation}\label{E:7}
(\alpha,\beta)={\tiny \begin{cases}
( 55, 55), &\quad \text{if }\quad  \chi(u^{11})=\chi(3a);\\
(61, 52), &\quad \text{if }\quad  \chi(u^{11})=\chi(3b); \\
(49, 58),  &\quad \text{if }\quad  \chi(u^{11})=2\chi(3a)-\chi(3b);\\
(43, 61),&\quad \text{if }\quad  \chi(u^{11})=3\chi(3a)-2\chi(3b);\\
(67, 49),&\quad \text{if }\quad  \chi(u^{11})=-\chi(3a)+2\chi(3b);\\
(37, 64),&\quad \text{if }\quad  \chi(u^{11})=4\chi(3a)-3\chi(3b).\\
\end{cases}}
\end{equation}
By  (\ref{E:2}) we obtain the system of inequalities
\[
\begin{split}
\mu_{0}(u,\chi_{3},*) & = \textstyle \frac{1}{33} (20(3\nu_{3b} + \nu_{11a}) + \alpha) \geq 0;\\ 
\mu_{11}(u,\chi_{3},*) & = \textstyle \frac{1}{33} (- 10(3\nu_{3b}+\nu_{11a}) + \beta) \geq 0,\\ 
\end{split}
\]
which  has no integer solutions in any of the six cases of (\ref{E:7}).

\noindent$\bullet$ Let $u$ be a unit of order $35$. By (\ref{E:1})
and Proposition \ref{P:4} we get $\nu_{5a}+\nu_{7a}+\nu_{7b}=1$.
Since $|u^{5}|=7$,  we need to
consider four cases for any character $\chi$ of $G$. They are  defined by part (v) of the Theorem.
By  (\ref{E:2}),  in all of the cases we get the system 
\[
\begin{split}
\mu_{0}(u,\chi_{2},*) & = \textstyle \frac{1}{35} (24(3 \nu_{5a} + 2 \nu_{7a} + 2 \nu_{7b}) + 47) \geq 0;\\ 
\mu_{7}(u,\chi_{2},*) & = \textstyle \frac{1}{35} (-6(3 \nu_{5a} + 2 \nu_{7a} + 2 \nu_{7b}) + 32) \geq 0,\\ 
\end{split}
\]
which has no integer solutions.

\noindent$\bullet$ Let $u$ be a unit of order $46$. By (\ref{E:1})
and Proposition \ref{P:4} we have that
$$
\nu_{2a}+\nu_{2b}+\nu_{23a}+\nu_{23b}=1.
$$
Put \quad $\alpha = \tiny { \begin{cases}
-3, \quad &\text{if }\quad \chi(u^{23})=\chi(2a); \\
1, \quad &\text{if }\quad \chi(u^{23})=\chi(2b);\\
9, \quad &\text{if }\quad \chi(u^{23})=-2\chi(2a)+3\chi(2b); \\
-7, \quad &\text{if }\quad \chi(u^{23})=2\chi(2a)-\chi(2b);\\
5, \quad &\text{if }\quad \chi(u^{23})=-\chi(2a)+2\chi(2b).\\
\end{cases} }
$

\noindent
According to (\ref{E:2}) we obtain that
\[
\begin{split}
\mu_0(u,\chi_2,3) &=-\mu_{23}(u,\chi_2,3)=\\
&=\textstyle\frac{1}{46} ( 22(
6\nu_{2a} -2\nu_{2b}-\nu_{23a}-\nu_{23b})+\alpha)= 0,\\
\end{split}
\]
which is impossible.

Now let \; $\chi(u^{23})=3\chi(2a)-2\chi(2b)$. \;
Put \; $t_1 = 3 \nu_{2a} - 5 \nu_{2b} +  \nu_{23a} +  \nu_{23b}$, \;
then by (\ref{E:2}) we obtain the system of inequalities
\[
\begin{split}
\mu_{23}(u,\chi_{3},*) & = \textstyle \frac{1}{46} (22 t_1 + 42) \geq 0;\quad 
\mu_{0}(u,\chi_{3},*)    = \textstyle \frac{1}{46} (-22 t_1 + 4) \geq 0, \\
\end{split}
\]
which has no solution for $t_1$.

\noindent$\bullet$ Let $u$ be a unit of order $55$. By (\ref{E:1})
and Proposition \ref{P:4} we have that  $\nu_{5a}+\nu_{11a}=1$.
By  (\ref{E:2}) we obtain the system of inequalities
\[
\begin{split}
\mu_{0}(u,\chi_{2},*) & = \textstyle \frac{1}{55} ( 40 (3 \nu_{5a} + \nu_{11a} ) + 45) \geq 0;\\ 
\mu_{11}(u,\chi_{2},*) & = \textstyle \frac{1}{55} (-10 ( 3 \nu_{5a} + \nu_{11a} ) + 30) \geq 0;\\ 
\mu_{1}(u,\chi_{2},*) & = \textstyle \frac{1}{55} (3 \nu_{5a} + \nu_{11a} + 19) \geq 0.\\ 
\end{split}
\]
It easy  to check that  last system of inequalities has no integral solution.

\noindent$\bullet$ Let $u$ be a unit of order $69$. By (\ref{E:1})
and Proposition \ref{P:4} we have that
$$
\nu_{3a}+\nu_{3b}+\nu_{23a}+\nu_{23b}=1.
$$
Since $|u^{23}|=3$ and by part (iv) of the Theorem we have six
cases for units of order 3, and, furthermore, $\chi(u^{3})\in\{\chi(23a),\chi(23b)\}$, 
for any character $\chi$ of $G$ we need to consider $12$ cases. Put
\begin{equation}\
(\alpha,\beta)={\tiny \begin{cases}
(23,23), &\; \text{if }\;  \chi(u^{23})=\chi(3a);\\
(29,20), &\; \text{if }\;  \chi(u^{23})=\chi(3b); \\
(17, 26),  &\; \text{if }\;  \chi(u^{23})=2\chi(3a)-\chi(3b);\\
(11, 29),&\; \text{if }\;  \chi(u^{23})=3\chi(3a)-2\chi(3b);\\
(35, 17),&\; \text{if }\;  \chi(u^{23})=-\chi(3a)+2\chi(3b);\\
(5, 32),&\; \text{if }\;  \chi(u^{23})=4\chi(3a)-3\chi(3b).\\
\end{cases}}
\end{equation}
By  (\ref{E:2}) in all of the $12$ cases we obtain the system
\[
\begin{split}
\mu_{0}(u,\chi_{3},*) & = \textstyle \frac{1}{69} (44(3 \nu_{3b} -  \nu_{23a} -  \nu_{23b})+ \alpha) \geq 0;\\ 
\mu_{23}(u,\chi_{3},*) & = \textstyle \frac{1}{69} (-22(3 \nu_{3b} -  \nu_{23a} -  \nu_{23b})+ \beta) \geq 0,\\ 
\end{split}
\]
which has no integral solutions.

\noindent$\bullet$ Let $u$ be a unit of order $77$. By (\ref{E:1})
and Proposition \ref{P:4} we have that
$$
\nu_{7a}+\nu_{7b}+\nu_{11a}=1.
$$
Since $|u^{11}|=7$,  we need to
consider four  cases for any character $\chi$ of $G$. They are  defined by part (v) of the Theorem.
By  (\ref{E:2}) we obtain the system of inequalities
\[
\begin{split}
\mu_{0}(u,\chi_{2},*) & = \textstyle \frac{1}{77} (60(2 \nu_{7a} + 2 \nu_{7b} +  \nu_{11a})+ 45) \geq 0;\\ 
\mu_{11}(u,\chi_{2},*) & = \textstyle \frac{1}{77} (-10(2 \nu_{7a} + 2 \nu_{7b} +  \nu_{11a})+31) \geq 0,\\ 
\end{split}
\]
which  has no integral solutions.

\noindent$\bullet$ Let $u$ be a unit of order $115$. By (\ref{E:1})
and Proposition \ref{P:4} we have that
$$
\nu_{5a}+\nu_{23a}+\nu_{23b}=1.
$$
Since $|u^{5}|=23$ and $\chi(u^{5})\in\{\chi(23a),\chi(23b)\}$,  we need to
consider two cases for any character $\chi$ of $G$. In both cases by  (\ref{E:2}) we get the system of inequalities
\[
\begin{split}
\mu_{0}(u,\chi_{2},*) & = \textstyle \frac{1}{115} (264 \nu_{5a} +
35) \geq 0;\qquad 
\mu_{23}(u,\chi_{2},*) = \textstyle \frac{1}{115} (-66 \nu_{5a} +
20) \geq 0, \\ 
\end{split}
\]
which  has no integral solution. The proof is done.

\noindent$\bullet$ Let $u$ be a unit of order $161$. By (\ref{E:1})
and Proposition \ref{P:4} we have that
$$
\nu_{7a}+\nu_{7b}+\nu_{23a}+\nu_{23b}=1.
$$
Since $|u^{23}|=7$ and $\chi(u^{7})\in\{\chi(23a),\chi(23b)\}$, for any character $\chi$ of $G$ we need to
consider eight   cases, defined by part (v) of the Theorem.
By  (\ref{E:2}) in all eight cases we obtain the system of inequalities
\[
\begin{split}
\mu_{0}(u,\chi_{2},*) & = \textstyle \frac{1}{161} (264 (\nu_{7a} + \nu_{7b} ) + 35) \geq 0; \\
\mu_{23}(u,\chi_{2},*) & = \textstyle \frac{1}{161} (-44( \nu_{7a} + \nu_{7b} ) + 21) \geq 0, \\
\end{split}
\]
which  has no integral solution.

\noindent$\bullet$ Let $u$ be a unit of order $253$. By (\ref{E:1})
and Proposition \ref{P:4} we have that
$$
 \nu_{11a}+\nu_{23a}+\nu_{23b}=1.
$$
Since  $\chi(u^{11})\in\{\chi(23a),\chi(23b)\}$,  we
consider two  cases for any character $\chi$ of $G$.
Put \quad $t_1 = 11 \nu_{23a} - 12 \nu_{23b}$ \quad and \quad 
$\alpha=\tiny{\begin{cases} 23 \quad &\text{if}\quad \chi(u^{11})=\chi(23a);\\
0 \quad &\text{if}\quad \chi(u^{11})=\chi(23b).\\
\end{cases}}$

\noindent
By (\ref{E:2}) in both  cases we obtain
\[
\begin{split}
\mu_{0}(u,\chi_{2},*) & = \textstyle \frac{1}{253} (220 \nu_{11a} + 33) \geq 0; \quad
\mu_{23}(u,\chi_{2},*)  = \textstyle \frac{1}{253} (-22 \nu_{11a} + 22) \geq 0;\\ 
\mu_{1}(u,\chi_{2},2) & = \textstyle \frac{1}{253} (t_1 + \alpha) \geq 0;\qquad \qquad 
\mu_{55}(u,\chi_{2},2)  = \textstyle \frac{1}{253} (- 10 t_1 + \alpha) \geq 0,\\ 
\end{split}
\]
so $\nu_{11a}=1$ and  $t_1=0$, so the solution is
$(\nu_{11a},\nu_{23a}, \nu_{23b})=(1,0,0)$. Now we compute that
$\mu_{1}(u,\chi_{2},*) = \textstyle \frac{1}{253} ( \nu_{11a} + 22)=\textstyle \frac{23}{253}$ 
is not an integer, thus, there is no solution in this case.

\bibliographystyle{plain}
\bibliography{Luthar_Passi_M24}

\vspace{30pt}
\centerline{\bf Appendix}
\vspace{5pt}
Possible partial augmentations $(\nu_{3a},\nu_{3b},\nu_{5a},\nu_{15a},\nu_{15b})$ for units of order 15:
$$
\tiny{
\begin{array}{llll}
( -3, 0, 5, -1, 0 ), &  ( -3, 0, 5, 0, -1 ), &  ( -2, -1, 5, -1, 0 ), &  ( -2, -1, 5, 0, -1 ), \\ 
( -2, 2, 5, -2, -2 ), &  ( -1, 1, 5, -3, -1 ), &  ( -1, 1, 5, -2, -2 ), &  ( -1, 1, 5, -1, -3 ), \\  
\bf{( 0, 0, 0, 0, 1 )}, & \bf{( 0, 0, 0, 1, 0 )}, &  ( 0, 3, 0, -1, -1 ), &  ( 1, -1, 0, 0, 1 ), \\
( 1, -1, 0, 1, 0 ), &  ( 1, 2, 0, -2, 0 ), &  ( 1, 2, 0, -1, -1 ), &  ( 1, 2, 0, 0, -2 ), \\  
( 2, 1, 0, -2, 0 ), &  ( 2, 1, 0, -1, -1 ), &  ( 2, 1, 0, 0, -2 ), &  ( 2, 4, 0, -3, -2 ), \\  
( 2, 4, 0, -2, -3 ), &  ( 3, 0, -5, 1, 2 ), &  ( 3, 0, -5, 2, 1 ), &  ( 3, 3, -5, 0, 0 ), \\
( 4, -1, -5, 1, 2 ), &  ( 4, -1, -5, 2, 1 ), &  ( 4, 2, -5, -1, 1 ), &  ( 4, 2, -5, 0, 0 ), \\  
( 4, 2, -5, 1, -1 ), &  ( 5, 1, -5, -1, 1 ), &  ( 5, 1, -5, 0, 0 ), &  ( 5, 1, -5, 1, -1 ), \\ 
( 5, 4, -5, -2, -1 ), &  ( 5, 4, -5, -1, -2 ). \\
\end{array}}
$$

\vspace{5pt}

Possible partial augmentations $(\nu_{3a},\nu_{3b}, \nu_{7a}, \nu_{7b}, \nu_{21a},\nu_{21b})$ 
for units of order 21:
$$
\tiny{
\begin{array}{lll}
( 0, 0, -3, 3, -1, 2 ), & ( 0, 0, -2, 2, 0, 1 ), & ( 0, 0, -1, 1, 0, 1 ),\\
\bf{( 0, 0, 0, 0, 0, 1 )}, & \bf{( 0, 0, 0, 0, 1, 0 )}, & ( 0, 0, 1, -1, 1, 0 ),\\
( 0, 0, 2, -2, 1, 0 ), & ( 0, 0, 2, -2, 2, -1 ), & ( 0, 0, 3, -3, 2, -1 ),\\
( 1, 2, -2, 2, -2, 0 ), & ( 1, 2, -1, 1, -2, 0 ), & ( 1, 2, -1, 1, -1, -1 ),\\
(1, 2, 0, 0, -1, -1 ), & ( 1, 2, 1, -1, -1, -1 ), & ( 1, 2, 1, -1, 0, -2 ),\\
( 1, 2, 2, -2, 0, -2 ), & ( 4, 2, -4, -3, 0, 2 ), & ( 4, 2, -4, -3, 1, 1 ),\\
( 4, 2, -3, -4, 1, 1 ), & ( 4, 2, -3, -4, 2, 0 ), & ( 0, 0, -2, 2, -1, 2 ).\\
\end{array}}
$$

\vspace{5pt}

\end{document}